\documentclass[a4paper,12pt]{article}
\usepackage{amsthm,amssymb}
\usepackage{graphicx}
\usepackage{subfigure}
\usepackage{hyperref}
\usepackage{subfigure}

\graphicspath{{images/}}%{misc/}}

\def\RR{\mathbb R}

\def\PP{\mathbb P}

\def\TT{\mathbb T}

\newtheorem{theorem}{Theorem}[section]

\theoremstyle{definition}

\theoremstyle{remark}

\begin{document}
\rmfamily
\title{Polyhedral Embeddings and Immersions of Many Triangulated $2$-Manifolds with Few Vertices}
\author{Ulrich Brehm, Undine Leopold}
\date{\today}
\maketitle
\begin{abstract}
\noindent This article presents an improvement and extension of the heuristic first presented by Hougardy et. al. \cite{hlz2010:srwtisf} for realizing triangulated orientable surfaces with few vertices by a simplex-wise linear embedding. The improvement consists in the applicability to non-orientable surfaces (simplex-wise linear immersions) as well as symmetric realizations. 
With the help of our algorithm, numerous - often symmetric - realizations of non-orientable surfaces with the minimal number or few vertices have been obtained for the first time. These examples include the  the Projective Plane with one or two handles and the Klein Bottle with one or two handles. \\[1ex]
 \textbf{Keywords:} Polyhedral Manifold, Polyhedral Surface, Triangulations, Geometric Realization, Symmetric Realization, Polyhedral Embedding, Polyhedral Immersion\\[0.5ex]
\textbf{Mathematics Subject Classification (2010):} Primary: 52B70, 52B55, Secondary: 52B10, 52B15, 52C40 

\end{abstract}

%%%%%%%%%%%%%%%%%%%
\section{Introduction}
%%%%%%%%%%%%%%%%%%%

Due to Steinitz' theorem \cite{s1922:pur, z1995:lop}, every polyhedral map on the sphere can be realized as the boundary complex of a convex $3$-polytope. In \cite{abem2007:htetmis}, it is proved that all triangulations on the torus can also be geometrically realized in $3$-space without self-intersection. For triangulations of orientable (closed) surfaces of higher genus, deciding the question of geometric realizability may present a challenging problem; the result by Bokowski and Guedes de Oliveira \cite{bgo2000:otgoom} proves that there are triangulations which are \emph{not} realizable as an embedded polyhedral surface. By Schewe's result \cite{s2010:nmvtossnuomass}, this holds true even for every genus greater than or equal to $5$. Moreover, a similar question can be posed for non-orientable triangulated surfaces, where geometric realization means immersion with flat, full-dimensional faces. In that context, an older result of the first author \cite{b1983:antms} shows that every non-orientable surface possesses triangulations which are not realizable in this way.

A naturally arising extremal problem is that of finding the minimal number of vertices sufficient for a geometric realization of a particular surface. The Heawood bound (see Section \ref{sec:background}) provides only a combinatorial lower bound.  
For a long time, computational search for realizations was too time-consuming, so examples have been constructed intuitively by hand, with Bokowski's \emph{rubber-band technique} \cite{b1989:agrwsidefdrm,b2008:ohmffros}, or with the help of a graphics software. See \cite{c1949:apwd,b1981:pmzevgd,b1987:amspog3w10v,bb1989:apog4wmnovams,b1990:htbmpmotbs,c1994:vmsiotkbits} for examples. 

In \cite{l2008:earrots}, Lutz first described a method for random realization of polyhedral embeddings of triangulated surfaces with few vertices. Random sets of integer coordinates were tested on their compatibility with the given triangulation. Furthermore, successful sets of coordinates were ``recycled'' or perturbed slightly and then tested for other triangulations. In a second paper, Hougardy, Lutz, and Zelke \cite{hlz2010:srwtisf} described a search heuristic on the integer lattice, decreasing an objective function measuring the length of the self-intersection to zero (signifying that an embedding has been reached).

Remarkably, both methods lead to success, with the search heuristic proving to be clearly superior over random realization. Intrigued, we adapted the algorithm presented in \cite{hlz2010:srwtisf} to deal with two new situations: immersed triangulations and symmetric realizations (for small symmetry groups). As input we used the complete lists of triangulations from Lutz's \emph{Manifold Page} \cite{lutz:manifoldpage}.

With a modified objective function we have been able to produce surprisingly many results, for example for triangulated projective planes and Klein bottles with $9$ vertices. Also, for the first time, polyhedral realizations of triangulations of the non-orientable surfaces of genus $3$, $4$, $5$, and $6$, with $9$, $10$, $10$, and $10$ vertices, respectively, have been obtained. As a consequence, the extremal problem of determining the minimal number of vertices needed for a geometric realization has been settled for the genera $3$, $5$, and $6$. Note that the minimality for the non-orientable surfaces of genus $3$ and $6$ is due to the Heawood bound \cite{h1890:mct} (see also Section \ref{sec:background}).  However, for the non-orientable surface of genus $5$, the second author proved that $9$ vertices do not suffice for a realization, and thus $10$ is minimal \cite{l2009:peaiot2m,l2016:nerfpiots}.

To the best of our knowledge, the only vertex-minimal realizations of non-orientable surfaces prior to our investigation were obtained by explicit construction; Brehm \cite{b1990:htbmpmotbs} constructed examples for the projective plane with $9$ vertices (and showed that fewer vertices do not suffice), and Cervone \cite{c1994:vmsiotkbits} found some examples for the Klein bottle using $9$ vertices (and proved that $8$, the combinatorial minimum, do not suffice for a geometric realization). 

Furthermore, with our adapted search heuristic symmetric realizations have been obtained in numerous cases. As we had hoped, imposing suitable symmetry conditions (compatible with subsets of the Euclidean integer lattice) and thus reducing the number of parameters sped up the search in many cases. For example, we found $35$ symmetric polyhedral realizations of triangulations of the orientable surface $M_5$ of genus $5$. This more than doubles the number of (un-symmetric) vertex-minimal examples of \cite{hlz2010:srwtisf}. Moreover, we also obtained symmetric polyhedral realizations (immersions) with $10$ vertices of triangulations of the \emph{non-orientable} surfaces $N_4$ of genus $4$, $N_5$ of genus $5$, and $N_6$ of genus $6$. Our results, which are described in detail in Section \ref{sec:nonorient-res} and Section \ref{sec:orient-res}, are summarized briefly in Table \ref{tab1}.

\begin{table}[htbp]
\begin{tabular}{c|c|c|c}
surface & $|V|$ (min) & realizable (total) & symmetry (no.)\\\hline
$N_1$ & $< 9 (6)$ & $0(20)$\cite{b1990:htbmpmotbs} & --\\
$N_1$ & $9(6)$\cite{b1990:htbmpmotbs} & $ \mathbf{\geq 50} (134)$ & $C_3\cite{b1990:htbmpmotbs} \mathbf{(\geq 3)}$\\
$N_2$ & $8(8)$ & $0(6)$\cite{c1994:vmsiotkbits} & --\\
$N_2$ & $9(8)$\cite{c1994:vmsiotkbits} & $\mathbf{\geq 90}(187)$ & $C_2 \cite{c1994:vmsiotkbits} (\mathbf{\geq 11)}$\\
$N_2$ & $10(8)$ & $\mathbf{\geq 76} (4462)$ & $\mathbf{C_2 (\geq 76)}$, incl. reflective\\
$N_3$ & $\mathbf{9}(9)$ & $ \mathbf{\geq 4}(133)$ & ?\\
$N_4$ & $9(9)$ & $?(37)$ & ?\\
$N_4$ & $\mathbf{10}(9)$ & $ \mathbf{\geq 112}(13657)$ & $\mathbf{C_2(\geq 112)}$\\
$N_5$ & $9(9)$ & $0(2)$\cite{l2016:nerfpiots} &-- \\
$N_5$ & $\mathbf{10}(9)$ & $ \mathbf{\geq 10}(7050)$ & $\mathbf{C_3(\geq 10)}$\\
$N_6$ & $\mathbf{10}(10)$ & $ \mathbf{\geq 1}(1022)$ & $\mathbf{C_2(\geq 1)}$\\\hline
$M_2$ & $10(10)$\cite{b1981:pmzevgd} & $865(865)$ \cite{l2008:earrots}, & $C_2\cite{b1981:pmzevgd}\mathbf{(\geq 46),C_3(\geq 7),}$\\
&&\cite{b2008:ohmffros}&$\mathbf{C_4(\geq 2), D_2(\geq 2)}$\\
$M_3$ & $10(10)$\cite{b1981:pmzevgd} & $20(20)$\cite{hlz2010:srwtisf}  & $C_2\cite{b1981:pmzevgd}(\mathbf{\geq 3}), \mathbf{C_3(\geq 2}),$\\
&&&$D_2(\geq 1)$\cite{b1987:amspog3w10v}\\
$M_4$ & $11(11)$ \cite{bb1989:apog4wmnovams}& $821(821)$\cite{hlz2010:srwtisf}& $C_2\cite{bb1989:apog4wmnovams}(\mathbf{\geq 28})$\\
$M_5$ & $12(12)$\cite{hlz2010:srwtisf} &$\mathbf{\geq 35}(751593)$ & $\mathbf{C_2(\geq 35)}$\\
$M_6$ & $12(12)$&$0(59)$\cite{s2010:nmvtossnuomass}  & --
\end{tabular}
\caption{Summary of results. New results are \textbf{boldfaced}. I.e., if $|V|$ is boldfaced, no results with this number of vertices have been known previously. For details see Sections \ref{sec:nonorient-res} and \ref{sec:orient-res}.}
\label{tab1}
\end{table}

The algorithm has also been successfully applied to obtain symmetric
realizations of the combinatorially regular $15$-vertex triangulation
$\{3,10\}_6$ of $M_6$. A separate paper \cite{bl2016:ascrpog6w15v} is devoted to this remarkable
result. $\{3,10\}_6$ belongs to the same family of regular triangulations as
Dyck's regular map $\{3,8\}_6$ (see \cite{cm1980:garfdg, b1989:agrwsidefdrm, b1987:msprodrm}).

%%%%%%%%%%%%%%%%%%%%%%%%%%%%%
\section{Background}\label{sec:background}
%%%%%%%%%%%%%%%%%%%%%%%%%%%%%

A \emph{$2$-manifold} $S$ is a Hausdorff space
which is locally homeomorphic to $\RR^2$ (sometimes this requirement is relaxed, allowing \emph{boundaries}). Compact $2$-manifolds are called \emph{closed surfaces}. Closed connected surfaces are classified abstractly as the $g$-fold connected sum of tori $\TT^2\#\ldots\#\TT^2$ in the orientable case (a sphere with $g$ handles), denoted $M_g$, or the $h$-fold connected sum of projective planes $\PP^2\#\ldots\#\PP^2$ in the non-orientable case (a sphere with $h$ cross-caps), denoted $N_h$. For the remainder of this article, the term \emph{surface} refers to a closed connected surface.

A \emph{triangulation of a surface $S$} is a two-dimensional simplicial complex $\Delta$ whose underlying topological space $|\Delta|$ is homeomorphic to $S$. Two triangles in $\Delta$ are \emph{adjacent} if they share an edge. A $k$-face and a $j$-face of the simplicial complex are \emph{incident} if one is a subset of the other, and a \emph{flag} is a maximal chain of mutually incident faces.  
Euler's formula states that 
\[\chi=f_0-f_1+f_2,\]
where $f_0$, $f_1$,$f_2$ are the numbers of vertices ($0$-faces), edges ($1$-faces), and faces ($2$-faces), respectively. For triangulations, we have $3f_2=2f_1$.
A triangulation of a surface with Euler characteristic $\chi$ requires at least 
\[f_0 \geq \left\lceil \frac{7+\sqrt{49 - 24 \chi}}{2}\right\rceil\]
vertices according to Heawood \cite{h1890:mct}. This bound was shown to be tight for all surfaces by Ringel \cite{r1955:wmdgnfimwdzk} and Jungerman and Ringel \cite{jr1980:mtoos}, with the only exceptions being the orientable surface $M_2$ and the non-orientable surfaces $N_2$ and $N_3$, all of which need one extra vertex.

It is a well-known fact from topology that orientable surfaces can be smoothly embedded into $\RR^{3}$, whereas non-orientable surfaces cannot. However, non-orientable surfaces can be smoothly immersed into $\RR^{3}$. 

The analogous polyhedral version of this question is much harder to decide. Let $\Delta$ be a triangulation of a surface $S=|\Delta|$, and let $V$ be the set of vertices. Each assignment $\psi:V\rightarrow\RR^3$ of coordinates to the set of vertices induces a simplexwise linear map $\phi_\psi\colon\left|\Delta\right|\rightarrow\RR^3$. By the terms \emph{polyhedral realization} or \emph{geometric realization} of $\Delta$ we mean an injective map $\phi_\psi$, i.e., a simplex-wise linear embedding, if $S$ is orientable, and a \emph{locally} injective map $\phi_\psi$, i.e., a simplex-wise linear immersion, if $S$ is non-orientable. In both cases, local injectivity assures that $i$-dimensional faces of $\Delta$ are mapped to flat, $i$-dimensional simplices in $\RR^3$ (points, line segments, and triangles).  

Polyhedral realizations of triangulated surfaces are interesting objects from the geometric, topological and combinatorial point of view. Since smooth surfaces can easily be discretized when allowing many vertices, realizations of triangulations with few or the minimal number of vertices present a greater challenge, and are thus of particular interest.

%%%%%%%%%%%%%%%%%%%%%%%%%%%%
\section{An Adapted Search Heuristic}
%%%%%%%%%%%%%%%%%%%%%%%%%%%%

We will now present the basic structure of the adapted algorithm, reviewing the approach in \cite{hlz2010:srwtisf}. Throughout this section, $\phi_\psi$ always denotes a simplexwise linear map induced by a coordinate assignment on the vertices, $\psi\colon V\rightarrow \mathbb{Z}^3$. Note that the coordinates, for the purpose of computational treatment, are \emph{integral}. Whether we mean the vertices, edges or triangles in the triangulation or whether we mean their images in $\RR^3$ can be derived from context. 

Comprehensive text files detailing (up to isomorphism) all possible triangulations of a specific closed surface with few vertices are readily available on Lutz's \emph{Manifold Page} \cite{lutz:manifoldpage} together with references pointing out their origins. These combinatorial
data were the basis for our systematic search for realizations.

Hougardy, Lutz, and Zelke \cite{hlz2010:srwtisf} modified vertex coordinates such that $\phi_\psi$ eventually becomes an embedding, all the while maintaining general position for the vertices. The goal of our extended algorithm is similar; however, we want that $\phi_\psi$ becomes an embedding for orientable surfaces, and an immersion for non-orientable surfaces. Furthermore, we would like to incorporate symmetry constraints. For these purposes we adapted the objective function as well as the general position requirement.

\begin{figure}[bth]
        \subfigure[Triangles not sharing vertices.]
        %{\label{fig:trianglesA}}
        {\includegraphics[width=0.65\linewidth]{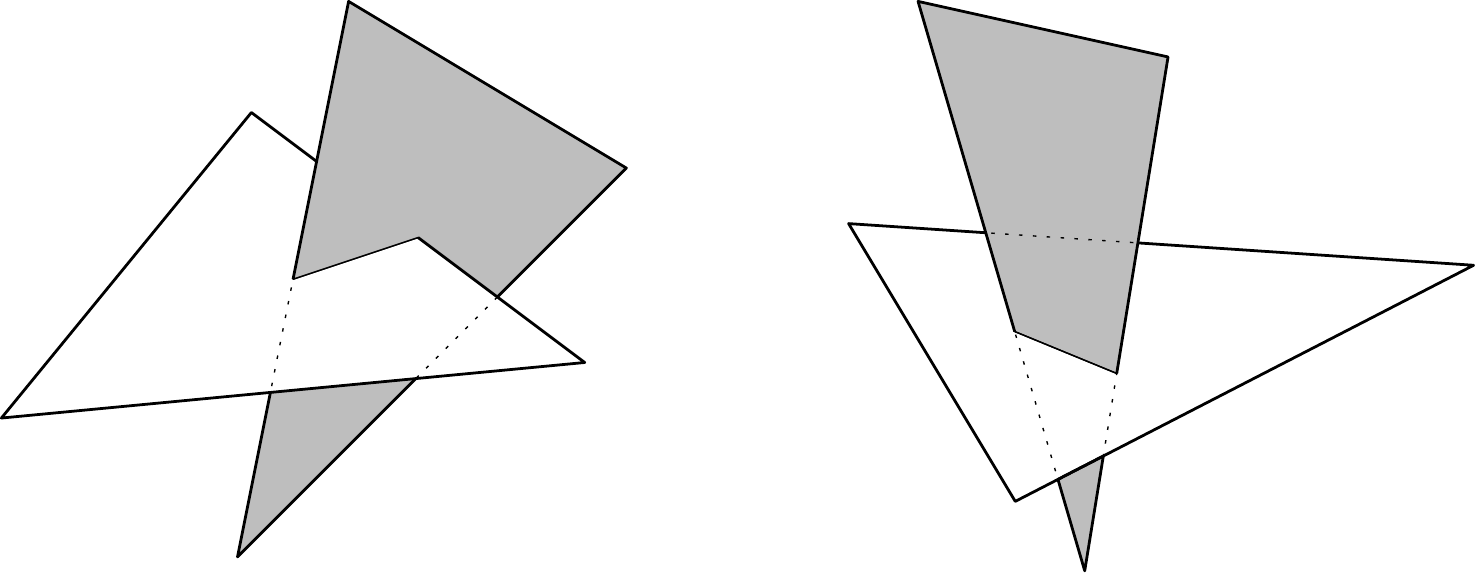}}\quad
        \subfigure[Triangles incident to the same vertex.]
        %{\label{fig:trianglesB}}
        {\includegraphics[width=0.25\linewidth]{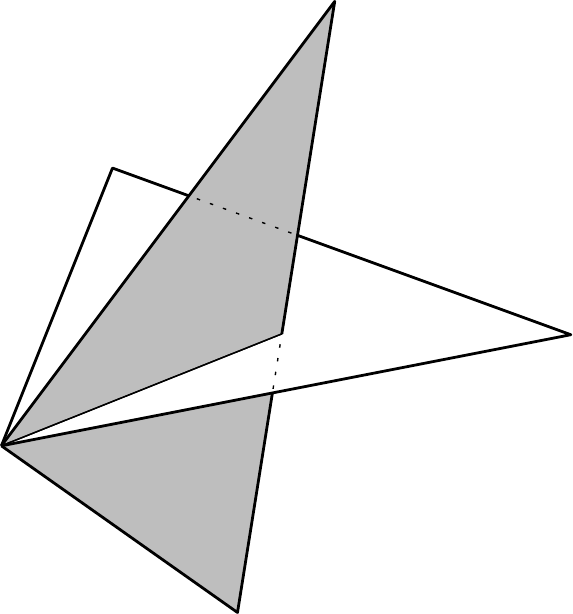}}
        \caption{Generic intersection of triangles in $3$-space.}
        \label{fig:triangles2}
\end{figure}

%----------------------------------------------------------------------------------------
\subsection{The Objective Function and Descent Step Heuristic}\label{subsec:iedge-emb}

When degenerate intersections are excluded, the intersection of two triangles is always a line segment, see Figure \ref{fig:triangles2}. Summing the lengths of the intersection segments for each pair of intersecting triangles, we obtain the value of the objective function used in \cite{hlz2010:srwtisf}. Note that it is not the absolute length of the self-intersection in $\RR^3$, as it cannot be excluded that $3$ triangles intersect in the same line segment. However, if we mark the intersection and then disassemble the surface into separate triangles, the objective function measures half of the total length of the markings. In what follows, we make use of this objective function when searching for \emph{embedded} realizations, and we modify it for immersions. 

The set of vertices of the triangulation is said to have \emph{admissible coordinates} if their coordinates are within a bounding box and fulfill some position requirement. A smaller bounding box is used for the initialization, and a bigger bounding box is used for all subsequent steps, allowing the surface to expand in the process. The specific position requirements are detailed in Subsection \ref{subsec:pos}.

A \emph{unit move} or \emph{unit coordinate move} is a change of value $+1$ or $-1$ in one of the coordinates of one of the vertices. An \emph{admissible step} or an \emph{admissible move} is a unit coordinate move which maintains admissible coordinates for the entire set of vertices. 
The descent step heuristic attempts to find a sequence of decreasing admissible moves.  More specifically, a descent step consists of the program cycling through the available unit moves in random order in search of admissible moves. After finding an admissible move, the objective function is re-computed. The first admissible move found is saved for the case that none of the admissible moves decrease the objective function. If a decreasing admissible move is found, it is taken and the program proceeds to the next descent step.

If no admissible move decreases the objective function, then a local minimum has been found. Hougardy et al. \cite{hlz2010:srwtisf} resolved this by using any admissible move and, in order to avoid endless undetected cycling between different local minima, a time limit was imposed up to which the objective function had to vanish, or new starting coordinates were set. By contrast, we only impose a maximum number of steps for each triangulation, and sometimes immediately re-initialize the coordinates when trapped in a local minimum. Specifically, if there is no decreasing move but there exists a saved admissible move, the program decides (pseudo-) randomly based on a given ratio on whether to take an admissible step or to restart with new coordinates. The step counter is not reset even though the coordinates are, so the program will eventually abort the search for each triangulation.

%--------------------------------------------------------------------------
\subsection{Position Requirements}\label{subsec:pos}
Numerical problems in calculating the objective function tend to arise from (near-) degenerate triangles or (near-) degenerate intersections of triangles.

General position of the vertices is the most obvious choice to assure that triangles in $\RR^3$ are full-dimensional and only intersect as in Figure \ref{fig:triangles2}. It can be verified by confirming affine independence of the coordinates for all quadruplets of vertices.

\begin{figure}[bth]
\begin{center}
        \subfigure[Vertex on a triangle.]
        {\includegraphics[width=0.4\linewidth]{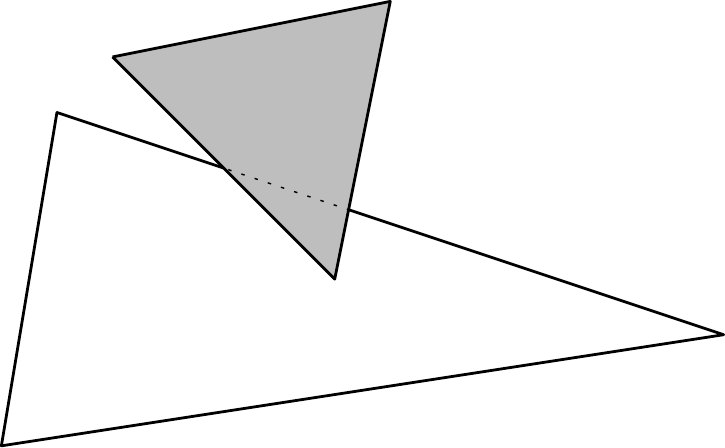}}\quad
        \subfigure[Edges intersecting in their relative interiors.]
        {\includegraphics[width=0.3\linewidth]{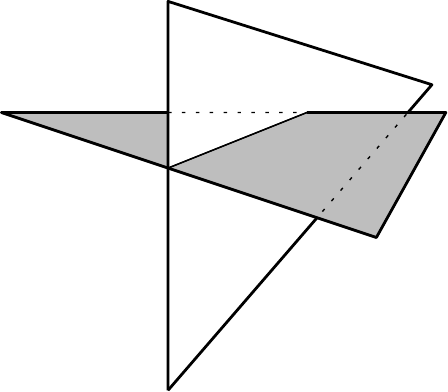}}
\end{center}
        \caption[Examples of Degenerate Intersections]{Examples of degenerate intersections.}
        \label{fig:degen-intersections}
\end{figure}

However, general position is more restrictive than necessary and would exclude some types of symmetry. For the purpose of symmetric realizations, a modified or relaxed version of the general position requirement for the vertices was used. In the context of this article, \emph{relaxed general position} of $\psi(V)$ can be characterized by $\psi$ being injective and preventing the degenerate intersections pictured in Figure \ref{fig:degen-intersections}. Specifically, if $v$ is a vertex, and $abc$ is a triangle whose vertices are different from $v$, then $\phi_\psi(v)=\psi(v)$ and $\phi_\psi(abc)$ are disjoint, and if $ab$ and $cd$ are edges whose vertices are pairwise disjoint, then $\phi_\psi(ab)$ and $\phi_\psi(cd)$ are disjoint. Note that this also prevents the intersection of triangles in the same plane.

We can still safely calculate the objective function, i.e., the self-intersection in $\RR^3$ is still composed of line segments, while allowing more flexibility in the position of the vertices. For arbitrary (non-symmetric) polyhedral realizations, there is no reason for using the relaxed version and thus, also in view of run-time considerations, we used the simpler condition of ordinary general position in the program.

As in \cite{hlz2010:srwtisf}, integer coordinates are required for the vertices and an outer bounding box for the coordinates is set. Integer coordinates assure that rescaling is not necessary as the object cannot shrink to arbitrarily small size, the calculations are simple, and numerical problems can be avoided. Last but not least, it is nice to have small integer coordinates for the vertices. Bounding boxes of a heuristically determined size are introduced for similar reasons, i.e., numerical safety, speed, limiting the number of (unsuccessful) search routes, and to obtain small integer coordinates if possible. 

%----------------------------------------------------------------------------
\subsection{Realizing Non-Orientable Closed Surfaces}

\subsubsection{A Modified Objective Function}\label{subsec:iedge-imm}
In order to handle immersions, the simplest modification of the objective function is to only compute the length of those intersection segments that are contradictory to an immersion. That is, we only sum the lengths of nontrivial intersections of triangles with a common vertex (see Figure \ref{fig:triangles2} (b)). Note that adjacent triangles may not intersect (except in their common edge) because we assure non-degeneracy. We use this altered objective function with the same descent step procedure as before. It is not a priori clear that this strategy should yield results, since a self-intersection in the simplicial neighborhood of a vertex constitutes a pinch point, and pinch points are usually hard to dissolve, as Cervone \cite{c1994:vmsiotkbits} pointed out. However, the varied results in Section \ref{sec:nonorient-res} confirm this as a suitable choice.

\subsubsection{A Simple Obstruction: Triple Points}\label{subsec:triple}

A result of Banchoff \cite{b1974:tpasois} stipulates that the number of \emph{triple points} of an immersed non-orientable closed surface \emph{in general position} has the same parity as its Euler characteristic. For a surface immersed in $3$-space, general position is the usual notion that distinct disks of the surface may intersect in $3$-space, but transversally to each other and at most three in the same point (triple point), in topologically the same way as coordinate planes intersect. Furthermore, it is required that there are only finitely many such triple points.

Our position requirements ensure transversality for any intersecting disks. It is always possible to perform a small perturbation of the vertices of a given geometric realization (possibly leaving the integer grid) such that another geometric realization is induced where no more than three disks intersect in a point, and three intersecting disks always intersect in the way coordinate planes do. There must be finitely many such triple points as there are finitely many faces (triangles). Thus a polyhedral realization of a non-orientable surface can be perturbed into another polyhedral realization such that the image of the surface is in general position, where it must have at least one triple point if its Euler characteristic is odd. 

Producing a triple point necessitates, at the very least, three disjoint triangles. Moreover, each pair $abc$, $def$ out of the three vertex-distinct triangles must be able to intersect in a non-degenerate way, thus exactly two out of their six combined edges must pierce the respective other triangle. Since intersections of the triangles incident to the same vertex are prohibited for immersions, the other triangle adjacent to such a piercing edge may not have its third vertex in $\{a,b,c,d,e,f\}$. Thus the triangles must satisfy the following necessary condition: at most $4$ of the $6$ triangles adjacent to $abc$ or $def$ (i.e., sharing an edge) have all three of their vertices in $\{a,b,c,d,e,f\}$.  
We automatically tested for this condition for triangulated surfaces with odd Euler characteristic (Cervone \cite{c1994:vmsiotkbits} used the term \emph{edge-cut-analysis} for testing similar conditions), and concluded polyhedral non-immersibility in many cases. Refined topological methods, incorporating the edge-cut analysis, are used in \cite{l2009:peaiot2m,l2016:nerfpiots} also for surfaces of even Euler characteristic.

%%%%%%%%%%%%%%%%%%%%%%%%%%%%%%%%%
\section{Realizations with Symmetries}\label{sec:sym}
%%%%%%%%%%%%%%%%%%%%%%%%%%%%%%%%%

Triangulations may possess nontrivial automorphisms, and it is natural to ask which automorphisms of a triangulation $\Delta$ can be realized as geometric symmetries. 
In this article, we consider geometric symmetries of the realizations $\phi_\psi:|\Delta|\rightarrow\mathbb{R}^3$, not merely symmetries of the images. That is, \emph{geometric symmetries} are isometries of $\mathbb{R}^3$ which fix $\phi_\psi(|\Delta|)$ and are induced by an automorphism of $|\Delta|$. 
 Our goal is thus the realization of symmetries within our program's framework, i.e.,  for vertices chosen in the integer lattice $\mathbb{Z}^3$. 

Subsets of the integer lattice $\mathbb{Z}^3$ in Euclidean $\RR^3$ are compatible with several symmetries, e.g., the reflection in a suitable plane (such as the $xy$-plane) as well as $2$-, $3$-, and $4$-fold rotation, point symmetry (central inversion), and some rotatory reflections. Higher symmetry can sometimes be achieved by combining these symmetries. However, numerical problems may arise and, for high symmetry, initialization of coordinates may not even be possible, so we kept only a small selection of symmetries to work with. Our program can handle reflections in a plane, rotations of order $2$, $3$, $4$, point symmetry (central inversion), rotatory reflections of orders $4$ and $6$ with angles $\frac{\pi}{2}$ and $\frac{\pi}{3}$, respectively, as well as the dihedral rotation group $D_2$ (three mutually orthogonal axes of order $2$). Note also that our program is targeted towards surfaces with few vertices, most of which have only few automorphisms to begin with; moreover, in any case, highly symmetric realizations seem to be extremely
      rare in our context. 

In order to handle symmetries in our program, a small alteration of the stepping procedure is necessary. After each coordinate move or choice of (pseudo-)random coordinates, and before re-evaluating the objective function, an adaptation function must restore the symmetry. As one vertex is moved, the positions of the other vertices in the same orbit are adapted accordingly. However, the overall structure of the descent step heuristic remains the same, as one can view the symmetry as part of the position requirement.

An automorphism of a triangulation is completely determined by its action on one of the $4f_1$ flags of the triangulation. During pre-processing of each triangulation, our program checks whether mapping a particular triangle in any way to any other triangle yields an automorphism, thus enumerating all of the automorphisms of the triangulation. 
 
In the main program, the non-trivial automorphisms as well as sets of generators of $D_2$ are tested successively, i.e. an instance of the search (descent step procedure) is started for each potential geometric realization of these symmetries. The size of ${\rm Aut}(\Delta)$ and the chosen generators are recorded in the results file.

The realizability of an automorphism as geometric symmetry depends, for example, on the vertices, edges, and faces fixed by the given automorphism. Some considerations are noted in the second author's diploma thesis \cite{l2009:peaiot2m}, and some necessary conditions have been incorporated into the program, pre-sorting automorphisms according to their potential geometric realizations. Thus, only the remaining candidates for
geometric symmetry are tested.  

Note that we are considering triangulations only and a plane square in $\mathbb{R}^3$ consisting of two coplanar adjacent triangles can never have $4$-fold symmetry according to our definition. Thus, the situation may be slightly different when we allow faces which have more than three vertices; however, for very few vertices, it seems unlikely that much more symmetric realizations could be found. For larger $|V|=f_0$, it is best to investigate beforehand which abstract symmetry groups have a chance of being realized geometrically in a given triangulation, thereby reducing the number of test cases. 

A well-known topological argument is that, generically, a line pierces an immersed surface in an even number of points. Since we allow singular points on rotation axes, an odd total number of fixed points by itself is not a criterion by which to exclude a symmetry (for immersions in the case that a rotation axis passes through a triple point). For related work, see \cite{t2014:tnomass, lt2015:esocsii3s}.

%%%%%%%%%%%%%%%%%%%%%%%%%%%%%%%%%%%%%%%%
\section{New Results for Non-Orientable Surfaces}\label{sec:nonorient-res}
%%%%%%%%%%%%%%%%%%%%%%%%%%%%%%%%%%%%%%

The results for triangulated non-orientable surfaces will be presented first. Their abundance confirms our modification of the objective function as suitable for accommodating non-orientable triangulated surfaces. The symmetric results indicate that symmetries can be realized heuristically as well. All results, both for the non-orientable and the orientable case, are available on the second author's web page \cite{leopold:home}. Unless otherwise stated, all results were obtained using at most $2000000$ search steps (descent steps) per triangulation. For the visualization we used JavaView \cite{javaview}.

%--------------------------------------------
\subsection{Projective Planes}

Considering that there has to be an odd number of triple points in a general position immersion, one easily concludes that rotations of even order cannot be geometric symmetries of a realization of $N_1$. 
We were able to generate many polyhedral realizations of triangulations of the projective plane with $9$ vertices, as well as a few with $3$-fold rotational symmetry (see also the examples found by the first author \cite{b1990:htbmpmotbs}). However, for almost as many triangulations their combinatorial structure prohibits triple points, rendering geometric realizations impossible.
 
\begin{theorem}
At least $50$ of the $134$ triangulations of the projective plane with $9$ vertices are polyhedrally realizable in $\RR^{3}$; at least $46$  are not, due to the infeasibility of a triple point. For at least $3$ triangulations there exist polyhedral realizations with a $3$-fold rotational symmetry.
\end{theorem}

For the search, the inner (initial) and outer bounding boxes had side lengths $40$ and $80$, respectively. 

Triangulations of the projective plane with one additional vertex, i.e. with $10$ vertices, were also processed, but only for symmetric realizations. With $60$ and $120$ for the side lengths of the inner and outer bounding box, we obtained $10$ realizations with $3$-fold rotational symmetry out of $1210$ triangulations (note that $923$ of the triangulations have trivial automorphism groups, see also \cite{l2008:earrots}). 

%---------------------------------
\subsection{Klein Bottles}

Cervone showed that the vertex-minimal triangulations (using $8$ vertices) are not geometrically realizable, and found some realizations of Klein bottles with $9$ vertices by hand, one of which possessed a $2$-fold rotational symmetry. We were able to considerably expand the list of Cervone's \cite{c1994:vmsiotkbits} examples. The bounding boxes had size $40$ (inner) and $80$ (outer). 
\begin{theorem}
At least $90$ of the $187$ triangulations of the Klein bottle with~$9$ vertices are polyhedrally realizable in $\RR^{3}$. At least $11$ triangulations possess polyhedral realizations with a $2$-fold rotational symmetry. 
\end{theorem}

\begin{figure}
\begin{center}
\includegraphics[width=.5\textwidth]{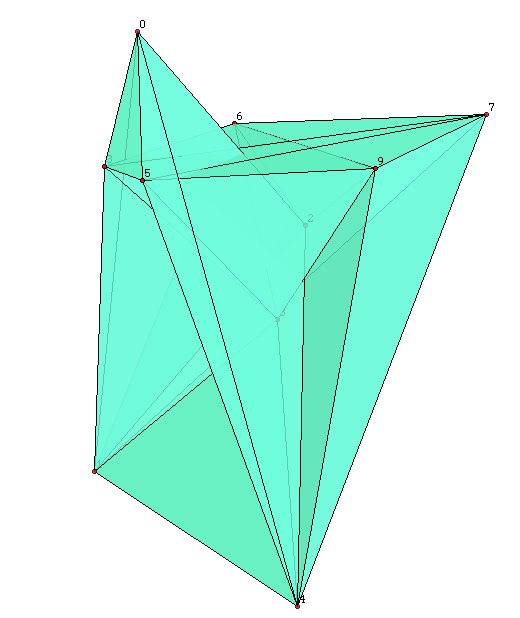}
\end{center}
\caption{A realization of $N_2$ with $10$ vertices and reflective symmetry.}
\label{fig:N2}
\end{figure}

In the hope of realizing further symmetries, we also tested $10$-vertex triangulations of the Klein bottle with larger bounding boxes (side length $60$ for the initial inner bounding box, side length $120$ for the outer bounding box) and the same default limit of $2000000$ on the search steps. This gave the following results.
\begin{theorem}
At least $76$ of the $4,462$ triangulations of the Klein bottle with $10$ vertices have polyhedral realizations in $\RR^{3}$ which realize an involutory symmetry. In most cases this was a rotational symmetry of order $2$, but in $9$ cases a reflective symmetry could be realized. 
\end{theorem}
Since a true reflective symmetry is difficult to realize for triangulations of surfaces with few vertices (the fact that there are only triangular faces being one obstruction), this is a rather remarkable result. We picture one realization in Figure \ref{fig:N2}.

%-------------------------------------------------------------------------
\subsection{Non-Orientable Surfaces of Higher Genus}

With the extended algorithm described in this article, realizations of the non-orientable surfaces of genus $3$, $4$, $5$, and $6$ with few or the minimal number of vertices could be obtained for the first time.

\medskip

\begin{figure}
\begin{center}
\includegraphics[width=.6\textwidth]{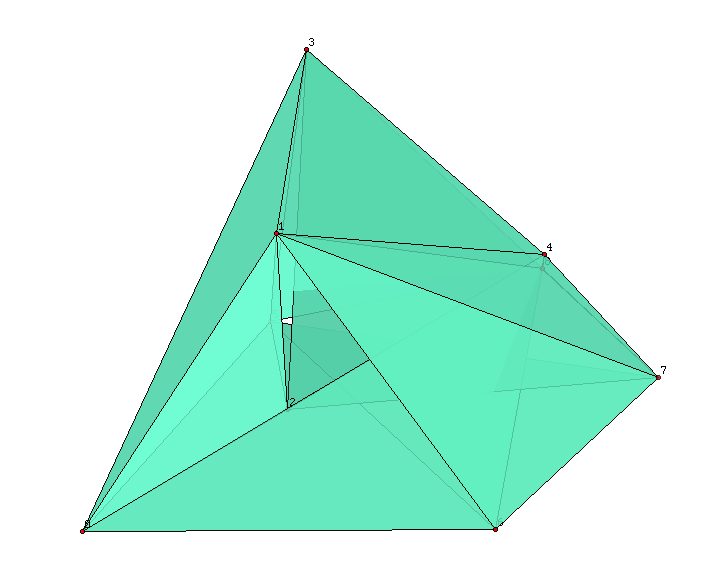}
\end{center}
\caption{A realization of $N_3$ with the minimal number of $9$ vertices.}
\label{fig:N3}
\end{figure}

First of all, the non-orientable surface $N_3$ of genus $3$ admits realizations with $9$ vertices, which is the minimal number of vertices needed for a triangulation due to the Heawood bound mentioned in Section \ref{sec:background}. One result of our search, using at most $8000000$ search steps and bounding boxes of sizes $60$ (initial) and $120$(outer), is shown in Figure \ref{fig:N3}. We did not find any symmetric realizations, but note that most ($106$ out of $133$, see also \cite{l2008:earrots}) of the triangulations possess no nontrivial automorphisms. Note also that not all triangulations may admit realizations.
\begin{theorem}
At least $4$ of the $133$ vertex-minimal $9$-vertex triangulations of the non-orientable surface of genus $3$ are polyhedrally realizable in $\RR^{3}$. At least $3$ of the $133$ triangulations are not polyhedrally realizable in $\RR^3$ due to the infeasibility of triple points.
\end{theorem}

\begin{figure}
\begin{center}
\includegraphics[width=.6\textwidth]{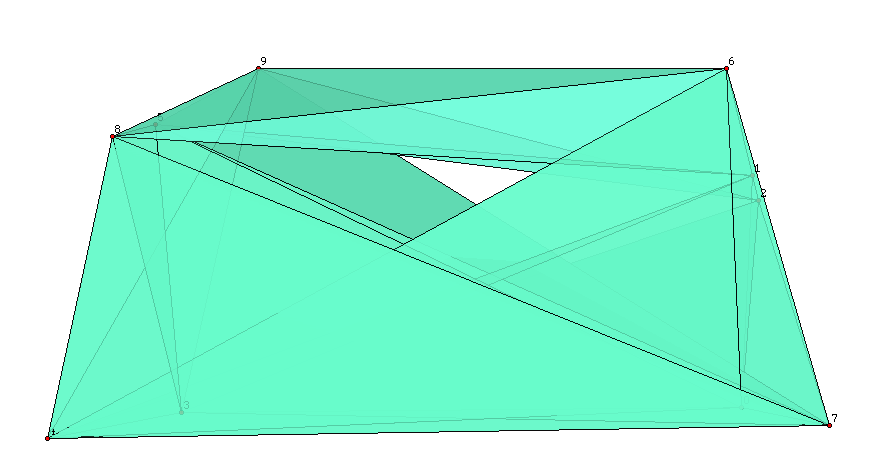}
\end{center}
\caption{A realization of $N_4$ with $10$ vertices and $2$-fold rotational symmetry.}
\label{fig:N4}
\end{figure}

For $N_4$, processing the vertex-minimal triangulations, i.e. triangulations with $9$ vertices, terminated without results. We used inner and outer bounding boxes of size $60$ and $120$, respectively. For the symmetric search we used the default maximum of $2000000$ search steps, for the asymmetric search we allowed even $10000000$. Still, no results were found and thus it remains an open problem whether $N_4$ can be realized with the minimal number of $9$ vertices.  However, by adding one vertex, we were able to obtain some, even symmetric, realizations for the non-orientable surface of genus $4$. We show one example in Figure \ref{fig:N4}.
\begin{theorem}
At least $112$ of the $13,657$ triangulations of the non-orientable surface of genus $4$ with $10$ vertices possess polyhedral realizations in $\RR^{3}$ with $2$-fold rotational symmetry.
\end{theorem}

The realizability of some $10$-vertex triangulations of the non-orientable surface of genus $5$ is one main result which we would like to highlight in this article. There exist two combinatorially distinct \emph{neighborly} triangulations (i.e., each pair of distinct vertices is connected with an edge) with $9$ vertices and none with fewer, due to the Heawood bound; however, the second author proved that these are not geometrically realizable \cite{l2009:peaiot2m, l2016:nerfpiots}, rendering $10$ minimal. The edge graph for a triangulation with $10$ vertices still has only $6$ missing edges. It is thus not only surprising that $10$ vertices suffice for a realization, but also, that symmetric realizations exist. 

\begin{theorem}
At least $10$ of the $7050$ triangulations of the non-orientable surface of genus $5$ with $10$ vertices admit polyhedral realizations in $\RR^3$. Moreover, all of the obtained realizations possess a $3$-fold rotational symmetry. 
\end{theorem}

\begin{figure}
\begin{center}
\includegraphics[width=.5\textwidth]{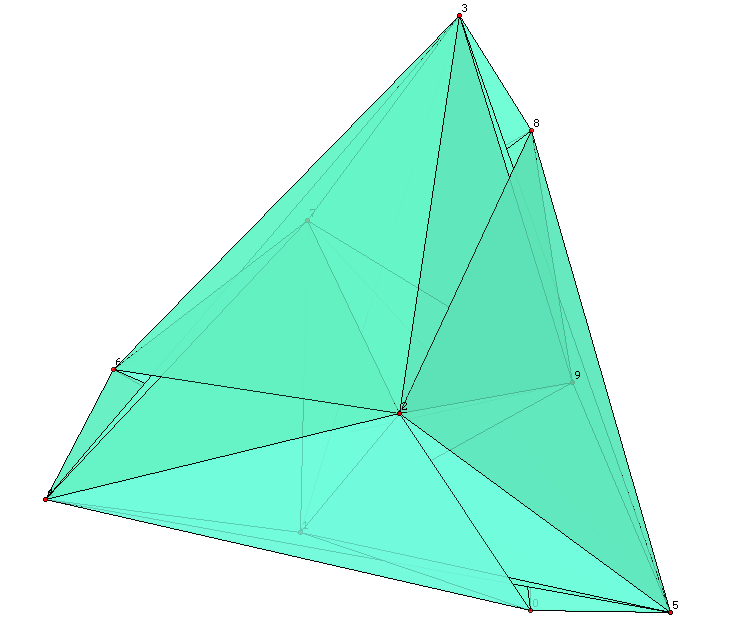}
\end{center}
\caption{A realization of $N_5$ with $10$ vertices and $3$-fold rotational symmetry.}
\label{fig:N5}
\end{figure}

All of the realizations took significantly less than $200000$ search steps, i.e. $\frac{1}{10}$ of the default limit. The side lengths of the inner (initial) and outer bounding boxes were $60$ and $120$, respectively. See Figure \ref{fig:N5} for an example.  

\medskip

\begin{figure}[htbp]
\begin{center}
\includegraphics[width=.7\textwidth]{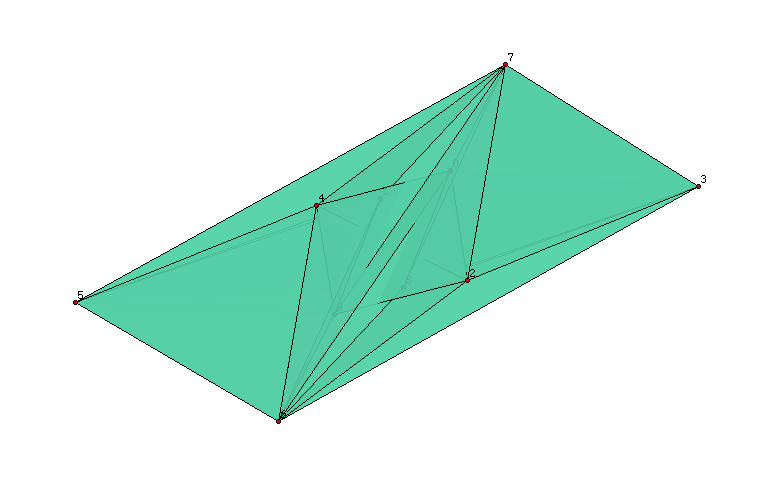}
\end{center}
\caption{A realization of $N_6$ with $10$ vertices and $2$-fold rotational symmetry.}
\label{fig:N6}
\end{figure}

Finally, we investigated vertex-minimal triangulations of $N_6$. One of the most remarkable results of our paper is that we
succeeded in finding one realization. This was achieved only
when searching for symmetric realizations. 
This realization is depicted in Figure \ref{fig:N6}.

\begin{theorem}
At least one of the $1022$  triangulations of the non-orientable surface of genus $6$ with the minimal number of $10$ vertices has a polyhedral realization in $\RR^3$ with $2$-fold rotational symmetry. 
\end{theorem}

%%%%%%%%%%%%%%%%%%%%%%%%%%%%%%%%%%
\section{Orientable Surfaces of Higher Genus}\label{sec:orient-res}
%%%%%%%%%%%%%%%%%%%%%%%%%%%%%%%%%%

The original results in \cite{l2008:earrots,hlz2010:srwtisf}, that all $865$ vertex-minimal triangulations of the orientable surface of genus $2$ are realizable in $\RR^3$, could be confirmed. 
We added to this by showing that several symmetries can be realized. We name the maximal symmetries with respect to each combinatorial type.

\begin{theorem}
At least $57$ of the $865$ vertex-minimal $10$-vertex triangulations of the orientable surface of genus $2$ have polyhedral realizations in $\RR^{3}$ which realize a non-trivial symmetry. Of these geometric realizations, $46$ possess $2$-fold rotational symmetry, $7$ possess $3$-fold rotational symmetry, $2$ realize a rotation-reflection of order $4$ and $2$ realize the dihedral group $D_2$ (Klein Four group generated by rotations). 
\end{theorem}

As parameters we used the default maximum of $2000000$ search steps and $40$ and $80$ for the size of the inner (initial) and outer bounding box, respectively. 

\medskip
The next candidates for symmetry were the realizations of the vertex-minimal triangulations of the orientable surface of genus $3$. Again, the realizability of all $20$ vertex-minimal triangulations of the orientable surface of genus $3$ \cite{hlz2010:srwtisf} could be confirmed. Moreover, our algorithm yielded the following symmetric versions, using an initial (inner) bounding box of side length of $40$ and allowing expansion up to side length $80$ for the outer bounding box. 

\begin{theorem}
At least $6$ of the $20$ vertex-minimal $10$-vertex triangulations of the orientable surface of genus $3$ have polyhedral realizations in $\RR^{3}$ which realize a symmetry. The maximal realized symmetries are the dihedral group $D_2$ ($3$ rotational axes of order $2$) for one triangulation, either a rotation of order $2$ or a rotation of order $3$ for one triangulation, a rotation of order $3$ for one other triangulation, and only a rotation of order $2$ for three more triangulations.
\end{theorem}
Brehm \cite{b1987:amspog3w10v} showed that the maximal order of the symmetry group is $4$ and gave an example with $D_2$ symmetry, and also gave an earlier example with $2$-fold rotational symmetry in  \cite{b1981:pmzevgd}. Moreover, the result was proved in the more general context of \emph{polyhedral realization of polyhedral maps}, i.e., where polygons with more vertices are allowed as faces (merging coplanar triangles and removing common edges). 

\medskip
Next, the vertex-minimal $11$-vertex triangulations of the orientable surface of genus $4$ were investigated. Bokowski and Brehm \cite{bb1989:apog4wmnovams} previously showed (again for polyhedral realizations of polyhedral maps) that the only possible nontrivial symmetry is a rotation of order $2$.
\begin{theorem}
At least $28$ of the $821$ vertex-minimal $11$-vertex triangulations of the orientable surface of genus $4$ have polyhedral realizations in $\RR^{3}$ with $2$-fold rotational symmetry.
\end{theorem}
This time, we used a higher maximum of $4000000$ search steps as well as inner and outer bounding boxes of side length $60$ and $120$, respectively.

\medskip
Finally, we were able to add symmetric versions to the results from \cite{hlz2010:srwtisf} concerning realizations of $12$-vertex triangulations of the orientable surface $M_5$ of genus $5$.

\begin{theorem}
At least $35$ of the $751,593$ vertex-minimal $12$-vertex triangulations of the orientable surface of genus $5$ admit polyhedral realizations in $\RR^{3}$ with $2$-fold rotational symmetry.
\end{theorem}

\begin{figure}[htbp]
\begin{center}
\includegraphics[width=.7\textwidth]{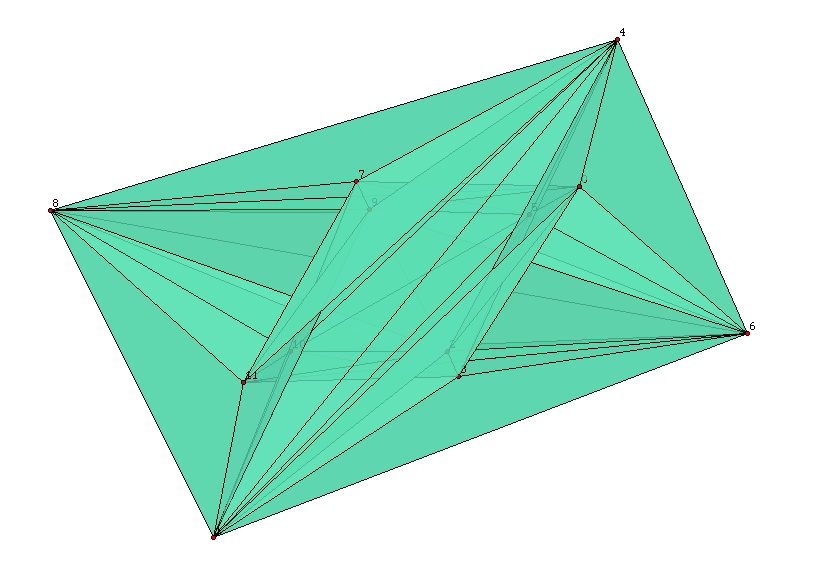}
\end{center}
\caption{A realization of $M_5$ with $12$ vertices and $2$-fold rotational symmetry.}
\label{fig:M5}
\end{figure}

Figure \ref{fig:M5} shows an example. The default of $2000000$ search steps was chosen as the maximal number of steps because of the high number of triangulations and the limited time frame. Several instances of the program ran independently on different portions of the file, and bounding boxes of size $60$ (inner) and $120$ (outer) were used. The requirement of geometric symmetry naturally limited the search to a subset of all $751,593$ triangulations. 

%%%%%%%%%%%%%%%%%%%%
\section{Final Remarks}
%%%%%%%%%%%%%%%%%%%

We have presented, implemented, and run a search heuristic adapted from the one given in \cite{hlz2010:srwtisf}.  It has proved to be suitable for finding simplex-wise linear realizations (i.e., embeddings or immersions) for triangulations of closed surfaces with few vertices. Our program is also capable of producing symmetric examples. Note that only small symmetries satisfying certain conditions have been tested, and the results are not exhaustive. Most of the results have already been obtained when the second author worked on her diploma thesis (see \cite{l2009:peaiot2m}) under supervision of the first author, and are published as a paper for the first time. 

The fact that the search terminates as soon as the intersection functional has reached zero may cause some pairs of edges to have very small distance. In the symmetric case, incorporating additional necessary conditions for geometric symmetries
may yield an improvement of the algorithm. 

The realizability problem for polyhedral surfaces is algorithmically decidable in theory \cite{bs1989:csg}, yet implementations of the procedure have only been successful in very specific cases (such as \cite{bb1987:anpog3w10v,s2010:nmvtossnuomass}). 
Recently, a refined approach has been taken by Firsching \cite{f2015:raifsssamp} for the successful realization of higher-dimensional complexes and spheres, and for proving non-realizability. 

The main open problem concerning the minimal number of vertices
for polyhedral realizations of surfaces of small genus is the case $N_4$,
where it still remains unclear if the minimal number of vertices for
a realization is $9$ or $10$.
In all other cases of orientable surfaces of genus up to $5$ and
non-orientable surfaces of genus up to $6$ we have found
vertex-minimal realizations. Moreover, in all these cases except for $N_3$
we have found symmetric realizations and except for $N_6$ we have
found several (or many) combinatorially different realizations.

The question whether \emph{all} triangulations of a certain surface are polyhedrally realizable is still open for $M_2$, $M_3$, and $M_4$ (and polyhedral embeddings). For non-orientable surfaces and polyedral immersions, Brehm's result \cite{b1983:antms} implies that there are non-realizable examples for any surface $N_h$.

%%%%%%%%%%%%%%%%%%%%
\bibliographystyle{amsalpha}
\bibliography{bibliography}

 \begin{flushleft}
Ulrich Brehm\\
Technische Universit\"at Dresden\\
Fachbereich Mathematik\\
Institut f\"ur Geometrie\\
D - 01062 Dresden\\
Germany\\[0.5ex]
ulrich.brehm@tu-dresden.de
\end{flushleft}
 
 \begin{flushleft}
Undine Leopold\\
Technische Universit\"at Chemnitz \\
Fakult\"at f\"ur Mathematik\\
D - 09107 Chemnitz\\
Germany\\[0.5ex]
undine.leopold@mathematik.tu-chemnitz.de
\end{flushleft} 

\end{document}